\documentclass[11pt,a4paper]{article}

\usepackage[utf8]{inputenc}
\usepackage{amsmath,amssymb,amsthm}
\usepackage{geometry}
\usepackage{hyperref}
\usepackage{cite}

\newtheorem{theorem}{Theorem}[section]
\newtheorem{lemma}[theorem]{Lemma}
\newtheorem{corollary}[theorem]{Corollary}

\theoremstyle{definition}
\newtheorem{definition}[theorem]{Definition}
\newtheorem{example}[theorem]{Example}

\title{\textbf{An Upper Bound on the Number of Distinct Composition Series in Finite Groups and Extremal Behavior in Group Isomorphism}}

\author{Abhijit Bhattacharjee\\
\small Email: \href{mailto:abhijit700121@gmail.com}{abhijit700121@gmail.com}}

\date{\today}

\begin{document}

\maketitle

\begin{abstract}
In this paper, we prove that among all finite groups of order $\le n$ (where $n \ge 4$ is a natural number), the number of distinct composition series is bounded above by $\prod_{i=1}^{\lfloor \log_2 n \rfloor} (2^i - 1)$, and equality holds if and only if $G$ is the elementary abelian 2-group of order $2^\alpha$, where $\alpha = \lfloor \log_2 n \rfloor$. This bound is a non-trivial upper bound of composition series and is best possible. Furthermore, we study the extremal asymptotic growth of this bound, proving it to be $n^{\frac{1}{2}\log_2 n + O(1)}$, and examine its direct implications for composition-series-based approaches to the group isomorphism problem.

\end{abstract}

\section{Introduction}

A composition series is a series of subgroups, each normal in the previous, such that the corresponding factor groups are simple. Any finite group has a composition series. The concept of a composition series in a group is due to Évariste Galois (1831) [Theorem 5.9, \cite{rotman}]. The famous Jordan--Hölder theorem, which states that the composition factors in a composition series are unique up to isomorphism, was proved in the nineteenth century \cite{jordan, holder}. Sometimes a group of small order has a huge number of distinct composition series. For example, an elementary abelian group of order $64$ has $615195$ distinct composition series. To find an upper bound for composition series in a finite group is a natural question. In \cite{hulpke}, there is an algorithm in GAP to find the distinct composition series of any group of finite order. The main aim of this paper is to provide an upper bound for the number of distinct composition series of any group of finite order using an elementary combinatorial approach, and to explore its implications in algorithmic group theory. All groups considered in this paper are of finite order.

There are five sections in this paper. In Section 2, we prove Theorem 2.1, which is the main preliminary theorem. In Section 3, we prove Theorem 3.7, which is a key counting result for finite abelian groups. In Section 4, we establish Theorem 4.3 by proving an underlying inequality in number theory, yielding the primary upper bound of the paper. Finally, in Section 5, we examine the asymptotic growth of this extremal bound and discuss its connection to worst-case complexity bounds in group isomorphism algorithms.

\section{Preliminary Results}

\begin{theorem}\label{thm2.1}
Among all finite groups of order $n$, the abelian group with elementary abelian Sylow subgroups has the highest number of distinct composition series.
\end{theorem}

Before proving Theorem \ref{thm2.1}, we recall some useful theorems in group theory.

\begin{theorem}[Correspondence Theorem, Theorem 2.28, \cite{rotman}]
Given $G \trianglerighteq N$, denote the family of subgroups of $G$ that contain $N$ by $\text{sub}(G; N)$ and the family of subgroups of $G/N$ by $\text{sub}(G/N)$. Then there is a bijection
\[
\phi : \text{sub}(G; N) \longrightarrow \text{sub}(G/N), \quad H \mapsto H/N
\]
which preserves all subgroup lattice and normality relationships.
\end{theorem}

\begin{theorem}[Exercise 14, Section 5.2, \cite{dummit}]
A finite abelian group is self-dual.
\end{theorem}

\begin{lemma}[Lemma 4.2.7, \cite{lauderdale}]
Suppose $G$ is a finite group and let $\mathcal{Q}(G) = \{p_1, p_2, \dots, p_r\}$ denote the distinct prime factors in the order of $G$. Let $P_i \in \text{Syl}_{p_i}(G)$. If $G = P_1 \times P_2 \times \dots \times P_r$, then
\[
|m(G)| = \sum_{i=1}^r |m(P_i)|,
\]
where $m(G)$ denotes the set of distinct maximal subgroups in $G$.
\end{lemma}

\begin{corollary}\label{cor2.5}
If $G$ is the abelian group with elementary abelian Sylow subgroups of order $n$, where $n = p_1^{\alpha_1} p_2^{\alpha_2} \dots p_r^{\alpha_r}$, then $m(G) = \sum_{i=1}^r \frac{p_i^{\alpha_i} - 1}{p_i - 1}$.
\end{corollary}

\begin{proof}
Since $G$ is a finite abelian group, it is self-dual and the direct product of its Sylow subgroups. Therefore, the number of maximal subgroups in $G$ is equal to the sum of all maximal subgroups in its Sylow subgroups. Now, an elementary abelian group is self-dual. So the number of maximal subgroups in an elementary abelian group is exactly equal to the number of minimal subgroups. Therefore, the number of subgroups of index $p_i$ is the number of distinct subgroups of order $p_i$ in an elementary abelian $p_i$-group of order $p_i^{\alpha_i}$, which is $\frac{p_i^{\alpha_i} - 1}{p_i - 1}$, and the result follows immediately.
\end{proof}

\begin{example}
If $G = \mathbb{Z}_2 \times \mathbb{Z}_2 \times \mathbb{Z}_2 \times \mathbb{Z}_2 \times \mathbb{Z}_3 \times \mathbb{Z}_3 \times \mathbb{Z}_5 \times \mathbb{Z}_5$, then $|G| = 3600$ and
\[
m(G) = \frac{2^4 - 1}{2 - 1} + \frac{3^2 - 1}{3 - 1} + \frac{5^2 - 1}{5 - 1} = 15 + 4 + 6 = 25.
\]
\end{example}

\begin{lemma}\label{lem2.7}
If $G = S_1 \times S_2$, where $S_1, S_2$ are two finite simple groups, and $G_1$ is a proper non-trivial normal subgroup of $G$, then $G_1$ is isomorphic to either $S_1$ or $S_2$.
\end{lemma}

\begin{proof}
We consider the projection map $p : G \to S_1$, which is onto. Therefore, the image of a normal subgroup under $p$ will be a normal subgroup in $S_1$. Since $S_1$ is simple, we have two cases:
\begin{description}
\item[Case 1:] If $p(G_1) = \{1\}$, then $G_1$ is a subgroup of $S_2$. Since it is a non-trivial normal subgroup of the simple group $S_2$, $G_1 = S_2$, i.e., $G_1 \cong S_2$.
\item[Case 2:] If $p(G_1) = S_1$, take the kernel $K = S_2 \cap G_1$ of the projection map. $K$ is normal in $S_2$, so either $K = \{1\}$ or $K = S_2$. If $K = \{1\}$, then $G_1 \cong S_1$. If $K = S_2$, then $G_1 \cong S_1 \times S_2$, a contradiction since $G_1$ is proper.
\end{description}
\end{proof}

\begin{corollary}\label{cor2.8}
If $G = S_1 \times S_2 \times \dots \times S_k$ is the direct product of $k$ non-abelian simple groups, then $G$ has $2^k$ normal subgroups, namely $N_1 \times N_2 \times \dots \times N_k$ where $N_i = \{1\}$ or $N_i = S_i$ for $1 \le i \le k$. Also, $G$ has only $k$ maximal normal subgroups, $M_1, M_2, \dots, M_k$, where $M_i = S_1 \times \dots \times S_{i-1} \times 1 \times S_{i+1} \times \dots \times S_k$ and $G/M_i \cong S_i$.
\end{corollary}

\begin{example}
Let $G = A_5 \times A_5$. Then $|G| = 3600$ and $G$ has $2^2 = 4$ normal subgroups, namely $1 \times 1, 1 \times A_5, A_5 \times 1, G$, and it has 2 maximal normal subgroups: $1 \times A_5$ and $A_5 \times 1$.
\end{example}

\begin{theorem}[Birkhoff, Theorem 8.6, \cite{burris}]
Every algebra $\mathcal{A}$ is isomorphic to a subdirect product of subdirectly irreducible algebras.
\end{theorem}

\begin{proof}[Proof of Theorem \ref{thm2.1}]
Let $G$ be the abelian group with elementary abelian Sylow subgroups such that $|G| = n$. By Lemma 2.4, the number of maximal subgroups in a finite nilpotent group of order $n$ is equal to the sum of distinct maximal subgroups in each Sylow subgroup. Among all finite $p$-groups of the same order, the elementary abelian group has the highest number of maximal subgroups. This proves that among all nilpotent groups of order $n$, $G$ has the highest number of maximal normal subgroups and hence the highest number of composition series.

Let $K$ be a solvable group of order $n$ with $m$ maximal normal subgroups $M_1, M_2, \dots, M_m$. The index of $M_i$ in $K$ is prime and $K/M_i \cong \mathbb{Z}_p$ for some prime $p$. Let $C'$ be the commutator subgroup of $K$; then $C' \unlhd M_i$. Let $M = \bigcap_{i=1}^m M_i$ be the intersection of all maximal normal subgroups of $K$. $M$ is a normal subgroup of $K$ such that $C' \unlhd M$, so $K/M$ is abelian. By the Correspondence Theorem, there exists a bijection $\phi : \text{sub}(K; M) \to \text{sub}(K/M)$ preserving lattice structure and normality. Thus $K/M$ is an abelian group with $m$ maximal normal subgroups. Since $G$ maximizes maximal normal subgroups among abelian groups of order $n$, $K$ cannot exceed $G$.

For a non-solvable, non-simple group, we consider the Jacobson radical $J(G)$ (the intersection of all maximal normal subgroups of $G$) and apply subdirect product decomposition from universal algebra. $J(G)$ is characteristic in $G$ and $J(G/J(G)) = 1$ \cite{guyot}. Thus $G/J(G)$ is a direct product of simple groups \cite{borovik}. If $J(G) = \{1\}$, $G$ is a direct product of simple groups. Otherwise, we obtain a minimal normal subgroup $P \le J(G)$ (the monolith) and quotient by $P$. Proceeding inductively, $G$ embeds into $G/N_1 \times \dots \times G/N_k$ where $\bigcap N_i = 1$ and each $N_i$ is maximal normal. If $H$ is a non-solvable, non-simple group of order $n$, $H$ is a subdirect product containing at least one non-abelian simple factor. By Corollary \ref{cor2.8}, a direct product of $k$ non-abelian simple groups has only $k$ maximal normal subgroups, completing the proof.
\end{proof}

\section{Counting the Composition Series in Finite Abelian Groups}

\begin{definition}
We define $\mathcal{C}_G$ as the set of all distinct composition series of the group $G$.
\end{definition}

\begin{theorem}\label{thm3.2}
Let $n \ge 2$ be a positive integer with prime factorization $n = p_1^{\alpha_1} p_2^{\alpha_2} \dots p_r^{\alpha_r}$. Then
\[
|\mathcal{C}_{\mathbb{Z}_n}| = \frac{(\sum_{i=1}^r \alpha_i)!}{\prod_{i=1}^r \alpha_i!}.
\]
\end{theorem}

\begin{proof}
For each positive divisor $d$ of $n$, there exists a unique subgroup of $\mathbb{Z}_n$ of order $d$. Let $S$ be the set of sequences of primes from $\{p_1, \dots, p_r\}$ of length $\alpha_1 + \dots + \alpha_r$, where $p_i$ occurs $\alpha_i$ times. Then $|S| = \frac{(\sum \alpha_i)!}{\prod \alpha_i!}$.

Define $f : S \to \mathcal{C}_{\mathbb{Z}_n}$ by
\[
f(\beta_1 \beta_2 \dots \beta_{\sum \alpha_i}) = \{e\} \unlhd \mathbb{Z}_{\beta_1} \unlhd \mathbb{Z}_{\beta_1 \beta_2} \unlhd \dots \unlhd \mathbb{Z}_{\beta_1 \dots \beta_{\sum \alpha_i}} = \mathbb{Z}_n.
\]
By construction, $f$ is injective. To show surjectivity, let $\{e\} = G_0 \unlhd G_1 \unlhd \dots \unlhd G_{\sum \alpha_i} = \mathbb{Z}_n$ be a composition series. Since every subgroup of $\mathbb{Z}_n$ is unique for a given order, $q_i = |G_i| / |G_{i-1}|$ is prime. Map $g : \mathcal{C}_{\mathbb{Z}_n} \to S$ by $g(\mathcal{C}_{\mathbb{Z}_n}) = q_1 q_2 \dots q_{\sum \alpha_i}$. Clearly $f^{-1} = g$, so $f$ is a bijection, and $|\mathcal{C}_{\mathbb{Z}_n}| = |S|$.
\end{proof}

\begin{example}
For $\mathbb{Z}_{360}$ with $360 = 2^3 \cdot 3^2 \cdot 5^1$, there are $\frac{(3+2+1)!}{3! 2! 1!} = 60$ distinct prime sequence permutations, each giving a unique composition series (e.g., $2 \cdot 3 \cdot 2 \cdot 2 \cdot 5 \cdot 3$ yields $\{e\} \unlhd \mathbb{Z}_2 \unlhd \mathbb{Z}_6 \unlhd \mathbb{Z}_{12} \unlhd \mathbb{Z}_{24} \unlhd \mathbb{Z}_{120} \unlhd \mathbb{Z}_{360}$).
\end{example}

\begin{theorem}\label{thm3.4}
Let $G$ be an abelian group of order $n = p_1^{\alpha_1} \dots p_r^{\alpha_r}$. Then
\[
|\mathcal{C}_G| = \prod_{i=1}^r t_i \cdot \frac{(\sum_{i=1}^r \alpha_i)!}{\prod_{i=1}^r \alpha_i!},
\]
where $t_i$ is the number of distinct composition series of the Sylow $p_i$-subgroup of $G$.
\end{theorem}

\begin{proof}
In a finite abelian group, composition factor orders correspond to prime factors. A composition series of $G$ induces composition series on each of its Sylow $p_i$-subgroups through canonical projection and subquotient relations. Because the Sylow $p_i$-subgroup has $t_i$ distinct composition series independently for each $i$, the total number of distinct composition series is obtained by multiplying the number of internal series assignments $\prod t_i$ by the total number of ways to interleave the prime factors, given by the multinomial coefficient $\frac{(\sum \alpha_i)!}{\prod \alpha_i!}$.
\end{proof}

\begin{theorem}\label{thm3.5}
In an elementary abelian group of order $p^n$, there are
\[
\prod_{k=1}^n \frac{p^k - 1}{p - 1} = \frac{\prod_{k=1}^n (p^k - 1)}{(p - 1)^n}
\]
distinct composition series.
\end{theorem}

\begin{proof}
This follows directly from Corollary \ref{cor2.5} by counting maximal subgroup chains recursively.
\end{proof}

\begin{theorem}\label{thm3.6}
Let $G$ be an abelian group with elementary abelian Sylow subgroups of order $n = p_1^{\alpha_1} \dots p_r^{\alpha_r}$. Then
\[
|\mathcal{C}_G| = \prod_{i=1}^r \left( \prod_{j=1}^{\alpha_i} \frac{p_i^j - 1}{p_i - 1} \right) \cdot \frac{(\sum_{i=1}^r \alpha_i)!}{\prod_{i=1}^r \alpha_i!}.
\]
\end{theorem}

\begin{proof}
This follows by combining Theorem \ref{thm3.4} and Theorem \ref{thm3.5}.
\end{proof}

\begin{theorem}\label{thm3.7}
Let $G$ be a finite group of order $n = p_1^{\alpha_1} \dots p_r^{\alpha_r}$. Then
\[
|\mathcal{C}_G| \le \prod_{i=1}^r \left( \prod_{j=1}^{\alpha_i} \frac{p_i^j - 1}{p_i - 1} \right) \cdot \frac{(\sum_{i=1}^r \alpha_i)!}{\prod_{i=1}^r \alpha_i!},
\]
and equality holds if and only if $G$ is an abelian group with elementary abelian Sylow subgroups.
\end{theorem}

\begin{proof}
This follows directly from Theorem \ref{thm2.1} and Theorem \ref{thm3.6}.
\end{proof}

\section{Bounding the Composition Series}

\begin{lemma}\label{lem4.1}
Let
\[
X = \left( \prod_{i=\alpha_1+1}^{\alpha_1+k} (2^i - 1) \right) \alpha_1! \alpha_r!, \quad \text{and} \quad Y = \left( \prod_{j=1}^{\alpha_r} \frac{p^j - 1}{p - 1} \right) (\alpha_1 + \alpha_r)!,
\]
where $p \ge 5$ is a prime, $\alpha_1 \ge 0$ is an integer, $\alpha_r \in \mathbb{N}$, such that $t = p^{\alpha_r}$ and $\lfloor \log_2 t \rfloor = k$. Assume $\alpha_r \ge 2$ whenever $p = 3$. Then $X/Y$ is a monotone increasing function of $\alpha_1$.
\end{lemma}

\begin{proof}
Define
\[
f(\alpha_1) = \frac{\prod_{i=1}^{\alpha_1+k}(2^i - 1)}{\prod_{i=1}^{\alpha_1}(2^i - 1) \prod_{j=1}^{\alpha_r}\left(\frac{p^j-1}{p-1}\right)} \cdot \frac{\alpha_1! \alpha_r!}{(\alpha_1 + \alpha_r)!}.
\]
Then
\[
\frac{f(\alpha_1 + 1)}{f(\alpha_1)} = \frac{2^{\alpha_1+k+1} - 1}{2^{\alpha_1+1} - 1} \cdot \frac{\alpha_1 + 1}{\alpha_1 + \alpha_r + 1}.
\]
To prove $f(\alpha_1 + 1) > f(\alpha_1)$, it suffices to show
\[
\frac{2^{\alpha_1+k+1} - 1}{2^{\alpha_1+1} - 1} > \frac{\alpha_1 + \alpha_r + 1}{\alpha_1 + 1}.
\]
Since $\frac{2^{\alpha_1+1+k} - 1}{2^{\alpha_1+1} - 1} > 2^k$ and $\frac{\alpha_1+\alpha_r+1}{\alpha_1+1} \le \alpha_r + 1$, it is sufficient to show $2^k > \alpha_r + 1 \iff k > \log_2(\alpha_r + 1)$.

Writing $\log_2 t = k + f$ where $0 < f < 1$, this inequality reduces to $\log_2\left(\frac{t}{\alpha_r + 1}\right) > f$. Since $f < 1$, it suffices to have $\frac{t}{\alpha_r + 1} > 2 \iff p^{\alpha_r} > 2\alpha_r + 2$, which holds for all $p \ge 5$ ($\alpha_r \ge 1$) and $p = 3$ ($\alpha_r \ge 2$) by induction. Thus $X/Y$ is monotone increasing in $\alpha_1$.
\end{proof}

\begin{theorem}\label{thm4.2}
Among all $p$-groups of order $\le n$ (where $n \ge 4$), the elementary abelian 2-group of order $2^\alpha$ (where $\alpha = \lfloor \log_2 n \rfloor$) has the highest number of composition series.
\end{theorem}

\begin{proof}
Among $p$-groups of a fixed order, the elementary abelian group maximizes composition series. Let $G_2$ be elementary abelian of order $2^\alpha$ with $\alpha = \lfloor \log_2 n \rfloor$, and $G_p$ be elementary abelian of order $p^\beta$ with $\beta = \lfloor \log_p n \rfloor$ ($p \ge 3$).
Since $2 \le p - 1$, we have
\begin{equation}\label{eq1}
2^{\lfloor \log_2 n \rfloor} - 1 > \frac{p^{\lfloor \log_p n \rfloor} - 1}{p - 1} \quad \text{for } n \ge 4.
\end{equation}
Taking products yields $|\mathcal{C}_{G_2}| = \prod_{i=1}^{\lfloor \log_2 n \rfloor}(2^i - 1) > \prod_{i=1}^{\lfloor \log_p n \rfloor}\frac{p^i - 1}{p - 1} = |\mathcal{C}_{G_p}|$. Thus $G_2$ maximizes the quantity among $p$-groups.
\end{proof}

\begin{theorem}\label{thm4.3}
Let $n \ge 4$ be a positive integer. Among all finite groups of order $\le n$, the elementary abelian group of order $2^\alpha$ has the highest number of composition series, where $\alpha = \lfloor \log_2 n \rfloor$. That is,
\[
|\mathcal{C}_G| \le \prod_{i=1}^{\lfloor \log_2 n \rfloor} (2^i - 1),
\]
and equality holds if and only if $G$ is the elementary abelian 2-group of order $2^\alpha$.
\end{theorem}

\begin{proof}
Let $H$ be an abelian group of order $m \le n$ with elementary abelian Sylow subgroups, $m = p_1^{\alpha_1} \dots p_r^{\alpha_r}$.

\begin{description}
\item[Step 1:] Let $p_r \neq 2$ with $p_r^{\alpha_r} = t$ and $k = \lfloor \log_2 t \rfloor$. Construct $H'$ of order $|H| \cdot 2^k / p_r^{\alpha_r} < |H|$.
\item[Step 2:] If $2 \nmid m$, then $k > \alpha_r$ and $\prod_{i=1}^k (2^i - 1) > \prod_{i=1}^{\alpha_r} \frac{p_r^i - 1}{p_r - 1}$, so $|\mathcal{C}_{H'}| > |\mathcal{C}_H|$.
\item[Step 3:] If $2 \mid m$, $H'$ has order $2^{\alpha_1 + k} p_2^{\alpha_2} \dots p_{r-1}^{\alpha_{r-1}}$. We compare $|\mathcal{C}_{H'}| / |\mathcal{C}_H|$. Using Lemma \ref{lem4.1}, the ratio is minimized at $\alpha_1 = 0$, reducing to the condition $\prod_{i=1}^k (2^i - 1) > \prod_{j=1}^{\alpha_r} \frac{p_r^j - 1}{p_r - 1}$, which holds by Theorem \ref{thm4.2}. Hence $|\mathcal{C}_{H'}| > |\mathcal{C}_H|$.
\item[Step 4:] Repeating this replacement step finitely many times reduces $H$ to an elementary abelian 2-group $G$ of order $2^\alpha \le n$. Finally, if $|H| = 2^{\alpha_1} \cdot 3^1$, replacing $\mathbb{Z}_3$ with $\mathbb{Z}_2$ increases the series count since $2^{\alpha_1+1} > \alpha_1 + 2$. Thus $G \cong (\mathbb{Z}_2)^\alpha$ achieves the unique maximum.
\end{description}
\end{proof}

\section{Extremal Structure and Algorithmic Complexity}

In this section, we analyze the asymptotic growth of the sharp upper bound established in Theorem \ref{thm4.3} and examine its implications for composition-series-based approaches to the group isomorphism problem \cite{rosenbaum}.

\subsection{Asymptotic Expansion of the Bound}

Taking the base-2 logarithm of the sharp upper bound $|\mathcal{C}_G| \le \prod_{i=1}^{\lfloor \log_2 n \rfloor} (2^i - 1)$, we obtain:
\[
\log_2 |\mathcal{C}_G| = \sum_{i=1}^{\lfloor \log_2 n \rfloor} \log_2(2^i - 1).
\]
Noting that $\log_2(2^i - 1) = i + \log_2(1 - 2^{-i}) = i + O(2^{-i})$, we evaluate the sum:
\[
\log_2 |\mathcal{C}_G| = \sum_{i=1}^{\lfloor \log_2 n \rfloor} i + O(1) = \frac{1}{2}\lfloor \log_2 n \rfloor^2 + O(\log n) = \frac{1}{2}(\log_2 n)^2 + O(\log n).
\]
Exponentiating both sides yields the precise asymptotic growth formula:
\[
|\mathcal{C}_G| = 2^{\frac{1}{2}(\log_2 n)^2 + O(\log n)} = n^{\frac{1}{2}\log_2 n + O(1)}.
\]

\subsection{Connection to Group Isomorphism Algorithms}

Algorithms that reduce the group isomorphism problem to composition-series isomorphism yield running time bounds that depend directly on the total number of composition series. In particular, standard composition-series enumeration frameworks bound the number of series by:
\[
|\mathcal{C}_G| \le n^{\frac{1}{2}\log_p n + O(1)},
\]
where $p$ is the smallest prime divisor of $|G|$.

Since $\log_2 n \ge \log_p n$ for any prime $p \ge 2$, the upper exponent is maximized precisely when $p = 2$. Theorem \ref{thm4.3} proves that this theoretical upper bound is achieved exclusively by elementary abelian 2-groups.

\subsection{Structural Interpretation}

This result provides a clear group-theoretic explanation for worst-case algorithmic complexity:
\begin{itemize}
    \item \textbf{Combinatorial Flexibility:} Elementary abelian $2$-groups $(\mathbb{Z}_2)^\alpha$ possess maximal subgroup lattice density, offering the largest possible choice of maximal normal subgroups at each step of a composition series chain.
    \item \textbf{Algorithmic Bottlenecks:} The exponent $\log_2 n$ appearing in general composition-series isomorphism bounds is not an artifact of proof techniques, but a direct consequence of the structural extremality of elementary abelian $2$-groups.
\end{itemize}

\section*{Acknowledgements}
The author utilized AI assistance for LaTeX formatting and refinement of English prose.

\end{document}